\documentclass[11pt]{amsart}
\usepackage[cp1251]{inputenc}
\usepackage[english]{babel}
\usepackage{amsmath}
\usepackage{amssymb}
\usepackage{amsfonts}
\numberwithin{equation}{section}
\begin{document}

\title[The Axially Symmetric Displacement Problem]{The Axially Symmetric Displacement Problem in the  Transversely Isotropic Elasticity}
\author{Yu. A. Bogan}
\address{Lavrent'ev Institute of Hydrodynamics, Lavrent'ev's prospect,  15, Novosibirsk, 630090, Russia}
\email{bogan@hydro.nsc.ru}

\maketitle

\begin{abstract}
In the assumption of hexagonal symmetry of an elastic material the axially symmetric displacement problem in a bounded axially symmetric solid with a Lyapunov boundary is reduced to a system of regular (Fredholm) integral equations. 
\end{abstract}

\section{Introduction}
\label{intro}
Here and throughout the paper the term 
 {\textquotedblleft axially symmetric\textquotedblright} refers to two properties: the first -- all relevant functions do not depend on the angular coordinate $ \phi $ of the cylindrical coordinate system  $ (r,\phi,z) $, and the second -- boundaries of all considered regions $ Q $  are  closed  surfaces of revolution about the $ z $-axis. 
Elastic problems involving axial symmetry are discussed at length in monographs \cite{Alex}, \cite{Pol}, \cite{Bel}. Elastic problems for transversely materials are discussed only in the book \cite{Alex}. The method of integral equations was applied  at first to problems of the axially symmetric elasticity in the dissertation of D. I. Sherman \cite{Sher}, which appeared in print in 1935; he solved there the axially symmetric isotropic traction problem. As the first step, he reduced it to singular integral equations, and as the second - to  Fredholm second kind equations. His computations are formidable in their complexity. He did not use the complex variable approach. The real variable approach was used later in \cite{Beck} by Becker. He remarked (\cite{Beck}, ch. 6, p. 125), that there are basically two approaches to axially symmetric formulations: the first is to derive from outset axially symmetric fundamental solutions based on ring loads as opposed to point loads and the second  --- to choose the three-dimensional point load solutions  and integrate them with respect to the angular coordinate. He remarked also, that both approaches lead to identical solutions. There is one obstacle, leading to serious mistakes when studying axially symmetric problems. It consists in the choice of the measure, entering into volume and, as a consequence, line integrals.  

 This paper is devoted to derivation of the correct representation of displacements in axially symmetric set-up for a transversely elastic material. The axially symmetric fundamental solutions based on ring loads were used for an isotropic material by Kermanidis in \cite{Ker}, and by Ishida and Ochiai in \cite{Ish1}, \cite{Ish2} for a transversely-isotropic material. There is a remarkable circumstance, which was left unnoticed by many researchers. It consists in the following: when one uses the traditional derivation of  the integral equation for a single elliptic second order equation, say, for a Dirichlet problem, based on the Green's formula and substitutes into it the fundamental solution, one arrives to the second kind Fredholm integral equation for an unknown density; nevertheless, the same approach, applied to an elliptic system, necessarily leads to a system of singular integral equations.  
We advances here the third approach, based on the simplicity of roots of the characteristic equation. This approach was successfully used earlier when studying two-dimensional elastic problems. A general solution of homogeneous equilibrium equations for a transversely elastic material depends on two different axially symmetric elliptic second order equations, differing only by real positive factors of second order derivatives with respect to $ z $.
The paper is organized as follows: 
\begin{enumerate}
\item As the first step, the axially symmetric half-plane displacement problem is solved in Section 3. Strictly speaking, this step (that is, solving the axially symmetric problem for a half-plane) is not obligatory; if one has in his disposal a correct notion of a Cauchy type axially symmetric integral, it is an easy task to derive the required system of integral equations;
\item the results of Section 3 are used later in Section 4 for deriving the complex-valued presentation of the previous solution and determination the correct form of the Cauchy type integral for a half-plane region;
\item Section 5 consists from two subsections: Subsection 5.1, where notion of the axially symmetric Cauchy-type integral is introduced and developed. Here is written the axially symmetric Cauchy-Riemann system of equations; in distinction from the usual standard complex variable theory, it consists from two different axially symmetric elliptic equations. The latter circumstance leads to difficulties in constructing the Cauchy-type integral. It is shown, that there are two candidates, which can be named the required {\textquotedblleft axially symmetric Cauchy type integral\textquotedblright}, but only one of them is the correct. It is essential, that the results of Section 4 permit us to choose the correct version. The results of Subsection 5.1 are used in Subsection 5.2; here Cauchy-type integrals  are written in regions $ Q_j, j=1,2 $ obtained from the region $ Q $ by suitable affine transformations. 
\item The main result of this paper is given in Section 6: a set of correct regular integral equations for the displacement problem. 
The author's aim in this paper is not to give one more ingenious representation of a solution of the axially symmetric elasticity (there are quite enough of them); but to give transparent formulas, leading to the correct real variable formulas for displacements. We hope, that this aim is achieved. The complex representation of a solution of the displacement problem is used only as an intermediate stage. 
The approach developed below for axially symmetric problems is in some respect similar to that of Sherman's, as that author also derived two real integral equations; in an axially symmetric bounded region with a sufficiently smooth boundary this approach always leads to a system of Fredholm second kind equations.
\end{enumerate}

The methods, used by Becker in \cite{Beck}, are quite sufficient for their numerical solution. 
It is rather strange, but there are no explicitly written exact solutions of standard axially symmetric boundary problems for a half-space, even for an isotropic material. The approach, developed in this paper, depends on the one hand, on the standard approach in modern partial differential equations, see \cite{Agmon}, and, on  the other hand, on simplicity of roots of the characteristic equation. The approach, developed in \cite{Agmon}, depends on the exact knowledge of Poisson kernels for a half-space for an arbitrary elliptic boundary value problem. The approach  of \cite{Agmon}  is modified below for an axially symmetric elastic problem.
In this paper is considered only the simplest case -- a displacement problem in the class of H\"{o}lder continuous boundary functions. The traction problem, which is somewhat more difficult, is left for a future paper. 
There are two interesting papers \cite{Ish1}, \cite{Ish2}, whose  authors write out the free-space Green's functions in \cite{Ish1}, and in \cite{Ish2} present a set of singular integral equations for the first boundary value problem. The integral equations, constructed below in this paper, do not coincide with the equations, presented in \cite{Ish2}. Recall also, that V. Chemeris  in \cite{Chem} solved the isotropic displacement problem by reducing it to a complex second kind Fredholm equation. 
 
\section{Preliminaries}
It is assumed, that in the cylindrical coordinate system $(r,\phi,z)$ the axis $z$ is taken orthogonal to the horizontal plane of isotropy, body forces are absent. Denote by  $\mathbf{R_{+}^2}= \{(r,z), z>0\}$ the upper half-plane. All relevant functions are extended to the negative semi-axis $ r<0 $ as even functions of $ r $: $ u(r,z)= u(-r, z), r >0 $. Recall, that an axial symmetry of a region does not change the generalized Hooke's law, since the stress-strain relations are indifferent to a form of a considered region.
The axially symmetric stress-strain relations are introduced as
\[ \sigma_{rr}=a_{11}\frac{\partial u_r}{\partial r }+a_{12}\frac{u_r}{r}
+a_{13}\frac{\partial u_z}{\partial z },\quad
\sigma_{zz}=a_{13}(\frac{\partial u_r}{\partial r}+\frac{u_r}{r})+a_{33}
\frac{\partial u_z }{\partial z },\]

\[ \sigma_{r \phi}= a_{12}\frac{\partial u_r}{\partial r }+
a_{11}\frac{u_r}{r}+a_{13}\frac{\partial u_z}{\partial z}, \quad
 \sigma_{rz}=a_{44}\Big( \frac{\partial u_r}{\partial z}+
\frac{\partial u_z}{\partial r }\Big),\]
where $ \sigma_{rr}, \sigma_{rz}, \sigma_{r \theta}, \sigma_{zz} $ are stresses, and $ u_r(r,z), u_z(r,z) $ are the components of the displacements vector with respect to the directions $ r $, $ z $  respectively; the component $ u_\phi $ in the direction $ \phi $ is equal to zero by the radial symmetry of  displacements. The matrix of elastic coefficients $ a_{ij}, i,j=1,2,6 $ is, as usual, positive definite.
It is well-known \cite{Eub}, that displacements can be written as linear combinations of first order derivatives of "quasi-harmonic" axially symmetric functions
\[ u_r(r,z)=\frac{\partial \varphi_1(r,z)}{\partial r}+
\frac{\partial \varphi_2(r,z)}{\partial r}, \quad
u_z(r,z)=k_1\frac{\partial \varphi_1(r,z)}{\partial z}+
k_2\frac{\partial \varphi_2(r,z)}{\partial z},\]
where a function $ \varphi_k(r,z), k=1,2$ is a solution of the equation
\begin{equation}
\frac{\partial^2 \varphi_k(r,z)}{\partial r^2}+\frac{1}{r}\frac
{\partial \varphi_k(r,z)}{\partial r }+\frac{1}{\lambda_k^2}\frac{\partial^2 \varphi_k(r,z)}
{\partial z^2 }=0. 
\end{equation}
Solutions of the previous equation are being called below quasi-harmonic functions. Here $ \lambda_1^2,\lambda_2^2$ are roots of the  equation
\[
a_{11}a_{44}+[a_{13}(a_{13}+2a_{44})-a_{11}a_{44}]\lambda^2 +a_{33}a_{44}\lambda^4=0. \]
It is assumed, that $ \lambda_1, \lambda_2 >0 $. Here
\begin{equation*}
k_1=\frac{a_{11}-a_{44}\lambda_1^2}{\lambda_1^2(a_{13}+a_{44})}, \quad
k_2=\frac{a_{11}-a_{44}\lambda_2^2}{\lambda_2^2(a_{13}+a_{44})}.
\end{equation*}
Then
\begin{equation*}
\sigma_{zz}=(a_{33} k_1-a_{13}\lambda_1^2 )\frac{\partial^2 \varphi_1(r,z) }
{\partial z^2}+(a_{33} k_2-a_{13}\lambda_2^2) \frac{\partial^2 \varphi_2(r,z) }
{\partial z^2},
\end{equation*}
\begin{equation*}
\sigma_{rz}=a_{44}(1+k_1)\frac{\partial^2 \varphi_1(r,z)}{{\partial r }{\partial z}}+
a_{44}(1+k_2)\frac{\partial^2 \varphi_2(r,z)}{{\partial r }{\partial z}}.
\end{equation*}
It is necessary to remind here the definition and  some properties of regular integral equations of the second kind.
\subsection{On the Fredholm integral equation}
Consider the integral equation
\[ f(s) + \lambda \int\limits_a^b K(s,s_0)f(s_0) \,d\,s_0=g(s), \]
where $ \lambda, a, b $ are  real parameters, $ f(s), g(s), K(s,s_0) $ are real functions. The function $ K(s,s_0) $ is defined in the plane $ (x,s) $ in the square $ a< s, s_0 < b $. According to the definition of S. G. Mikhlin \cite{Mikh}, the equation above is a Fredholm  equation  of second kind for the function $f(s)$, if $ g(s)$, and $K(s,s_0) $ are square integrable in the square $ a< s, s_0 < b $.
Recall, that Fredholm assumed continuity of the kernel $ K(s,s_0) $ in the same square. If a boundary of a plane region is a Lyapunov curve, then the kernel of the integral operator
\[ \frac{1}{\pi i}\int\limits_{\partial Q}\frac{f(s) \,d\,t}{t - z}, \]
where $f(s)$ is a real function, 
\[ t=t(s)= x_1(s)+ix_2(s), \quad \,d\,t=(x_1'(s)+i x_2'(s))\,d\,s,\quad z=x_1+i x_2, \]
is a kernel with a weak singularity, i.e., it can be written as a fraction
\[ K(s-s_0)=\dfrac{ a(s,s_0)}{|s-s_0|^\alpha}, 0< \alpha <1, \]
where $ a(s,s_0) $ is a bounded function. It can be proved, \cite{Mikh}, that if a kernel of an integral operator has a weak singularity, all iterated kernels, beginning  from some, are bounded. Hence, the equations with weakly singular kernels are  the Fredholm ones. On the other hand, the proof of uniqueness of the  axially symmetric elastic problem, as it is given in the book \cite{Nech}, uses only a piece-wise smoothness of a boundary; it means that the Fredholm alternative is true for regions more general than bounded regions with a Lyapunov boundary. 

The integral equations, written below for the axially symmetric elastic displacement problem, are the Fredholm ones for an axially-symmetric region with a Lyapunov boundary.

\section{The half-plane displacement problem}
We will give in this section the formal solution of the half-plane displacement problem, applying the Fourier-Hankel transform to quasi-harmonic functions. The explicit expressions of displacements, derived below, are sums of axially symmetric potentials. They will be used later to obtain the correct complex-valued representation of displacements in an arbitrary axially symmetric region. Recall, that the Hankel integral transform of order $ \nu $ is defined as
\[ h_{\nu} (f)=\int\limits_0^\infty f(x)J_{\nu}(xy)x \,d\,x, \]
where $J_\nu(y)$ is the first kind Bessel function of order $\nu$.  Hankel was the first to prove the inversion formula for it; recall \cite{Wats}, that the Hankel transform of a function $f(r)$ is valid at every point at which $ f(r) $ is continuous, provided that it is defined on $ (0, \infty) $, has a bounded variation in every finite subinterval in $ (0, \infty) $, and
\[ \int\limits_0^\infty |f(r)| \sqrt{r}\,d\,r \]
is finite. Then the inversion formula
\[ \int\limits_0^\infty u \,d\,u \int\limits_0^\infty f(x) J_\nu(ux)J_\nu(ur) x \,d\,x= f(x) \]
holds.
Recall also, that
\[ \int\limits_0^\infty J_\nu(kr) J_\nu(k'r)r \,d\,r =\dfrac{\delta(k-k')}{k}, \]
where $ \delta(k-k') $ is the one-dimensional $ \delta $-function.
Assign displacements  on the boundary of upper half-space ${\mathbf R}_+^2= \{(r,z), z >0 \}$ :
\begin{gather}
\Big(\frac{\partial \varphi_1(r, \lambda_1 z)}{\partial r}+
\frac{\partial \varphi_2(r, \lambda_2 z)}{\partial r}\Big)\vert_{z =+0}= f_1(r),\\
\Big(k_1\frac{\partial \varphi_1(r,\lambda_1 z)}{\partial z }+
k_2\frac{\partial \varphi_2(r,\lambda_2 z)}{\partial z }\Big)\vert_{z=+0}=f_2(r),
\end{gather}
where $ f_k(r), k=1,2 $  are continuous functions with bounded variation in any finite sub-interval of $ (0,\infty) $, which obey the inequalities
\[ \int\limits_0^\infty |f_k(r)| \sqrt{r}\,d\,r < \infty, \qquad k=1,2.\]
Apply the zero order Hankel transform to equations (2.1). Omitting some simple computations, we can write the functions $\varphi_k(r,z), k=1,2$ as
\begin{equation*}
\varphi_k(r,z) =\int\limits_{0}^{\infty} A_k(t)e^{-t\lambda_k z}t J_0(tr)\,d\,t, \qquad k=1,2.
\end{equation*}
Then the functions $ A_k(t),k=1,2$ are determined from
\begin{equation*}
-\sum\limits_{k= 1}^{2}\int\limits_{0}^{\infty}A_k(t)e^{-t\lambda_t z}t^2 J_1(tr)\,d\,t
\vert_{z=+0}=f_1(r),
\end{equation*}

\begin{equation*}
-\sum\limits_{k=1}^{2}\int\limits_{0}^{\infty}A_k(t)\lambda_ke^{-t\lambda_t z}t^2
J_0(tr)\,d\,t\vert_{z=+0}=f_2(r).
\end{equation*}
Put
\begin{equation*}
f_1(r)=\int\limits_{0}^{\infty}\hat{f}_1(t)t J_1(tr)\,d\,t, \quad
f_2(r)=\int\limits_{0}^{\infty}\hat{f}_2(t)t J_0(tr)\,d\,t.
\end{equation*}
Then $A_k(t), k=1,2$ are
\begin{equation*}
\delta A_1(t)t=k_2\lambda_2 \hat{f}_1(t)- \hat{f}_2(t), \quad
\delta A_2(t)t= \hat{f}_2(t)-k_1\lambda_1 \hat{f}_1(t),\quad
\delta=k_1\lambda_1-k_2\lambda_2,
\end{equation*}
and $\varphi_k(r,z), k=1,2$ are expressed as
\begin{equation*}
\varphi_r(r,z)=\frac{1}{\delta}\int\limits_{0}^{\infty}(k_2 \lambda_2
\hat{f}_1(t) -\hat{f}_2(t)) e^{-t\lambda_1 z} J_0(tr) \,d\,t,
\end{equation*}
\begin{equation*}
\varphi_z(r,z)=-\frac{1}{\delta}\int\limits_{0}^{\infty}(k_1 \lambda_1
\hat{f}_1(t) -\hat{f}_2(t)) e^{-t\lambda_2 z} J_0(tr) \,d\,t.
\end{equation*}
Finally, displacements are 
\begin{align}
u_r(r,z)= &-\frac{1}{\delta}\int\limits_{0}^{\infty}(k_2 \lambda_2 \hat{f}_1(t)
 -\hat{f}_2(t))J_1(tr)t e^{-t\lambda_1 z}\,d\,t +\\ \nonumber
&\frac{1}{\delta}\int\limits_{0}^{\infty}(k_1 \lambda_1 \hat{f}_1(t)-\hat{f}_2(t))
J_1(tr)t e^{-t\lambda_2 z}\,d\,t.
\end{align}
\begin{align}
u_z(r,z)= &-\frac{k_1}{\delta}\int\limits_{0}^{\infty}(k_2 \lambda_2 \hat{f}_1(t) -\hat{f}_2(t))
J_0(tr)\lambda_1 t e^{-t\lambda_1 z}\,d\,t +\\ \nonumber
&\frac{k_2}{\delta}\int\limits_{0}^{\infty}(k_1 \lambda_1 \hat{f}_1(t)-\hat{f}_2(t))
J_0(tr)t \lambda_2 e^{-t\lambda_2 z}\,d\,t.
\end{align}
Now, substitute the previous expressions for $ \hat{f}_k(t), k=1,2$ into (3.3), (3.4). As
\begin{equation*}
\hat{f}_1(t)=\int\limits_{0}^{\infty}a J_1(at)f_1(a) \,d\,a,
\quad \hat{f}_2(t)=\int\limits_{0}^{\infty}a J_0(at)f_2(a) \,d\,a,
\end{equation*}
the double integral
\begin{equation*}
A= \int\limits_{0}^{\infty}\int\limits_{0}^{\infty} t e^
{-t\lambda_k z}J_0(tr)\,d\,t aJ_0(at)f_2(a)\,d\,a
\end{equation*}
can be written as
\begin{equation*}
A= \int\limits_{0}^{\infty} \Lambda_k(r,z;a,0) f_2(a) a \,d\,a,
\end{equation*}
where
\begin{equation*}
\Lambda_k(r,z;a,0)= \int\limits_{0}^{\infty}J_0(at)J_0(tr) e^{-t\lambda_k z}t\,d\,t.
\end{equation*}
Now, $ \Lambda_k(r,z;a,0) $ is the derivative of the function 
\begin{equation}
g_k(r,z;a,0) = \int\limits_{0}^{\infty}J_0(at)J_0(tr)
 e^{-t\lambda_k z} \,d\,t,
\end{equation}
the axially symmetric fundamental solution of the equation (2.1) with pole $ (a,0) $.
It also can be written as
\[ g_k(r,z;a,0)=\dfrac{1}{\pi}\int\limits_0^\pi\dfrac{1}{\sqrt{r^2+a^2+\lambda_k^2 z^2-2ar\cos \alpha}} \,d\,\alpha.  \]
Therefore, $ \Lambda_k(r,z;a,0) $  is the axially symmetric Poisson kernel
\[ \Lambda_k(r,\lambda_k z;a,0)= \dfrac{z}{\pi}\int\limits_0^\pi \dfrac{1}{(r^2+a^2+\lambda_k^2 z^2-2a r\cos \alpha )^{3/2}} \,d\,\alpha \]
for the equation (2.1). The term \textquotedblleft Poisson kernel\textquotedblright  means, that if $ \varphi_k(r,\lambda_k z)\arrowvert_{z= +0}=h(r) $, then the solution of the Dirichlet problem for the quasi-harmonic equation (2.1) is represented as
\begin{equation*}
\varphi_k(r,\lambda_k z)= \int\limits_0^\infty \Lambda_k(r,z;a,0)h(a) a\,d\,a.
\end{equation*}
In particular,
\begin{equation*}
\Lambda_k(r,z;a,0)=-\frac{\partial }{\lambda_k \partial z} \int\limits_{0}^{\infty}J_0(at)J_0(tr)
 e^{-t\lambda_k z} \,d\,t =-\frac{\partial g_k(r,\lambda_k z;a,0)}{\lambda_k \partial z}, k=1,2.
\end{equation*}
Consider now the integrals
\[ \int\limits_0^\infty J_1(at)J_1(tr)t e^{-t\lambda_k z}\,d\,t\, k=1,2.\]
They can be written as follows:
\[ \frac{\partial}{\lambda_k \partial z} \int\limits_0^\infty J_1(at)J_1(tr) e^{-t\lambda_k z}\,d\,t\,=-\frac{\partial}{\lambda_k\partial z}G_k(r,\lambda_k z; a,0), k=1,2, \]
where
$ G_k(r, \lambda_k z;a,0) ,k=1,2 $ are the fundamental axially symmetric solutions of equations
\begin{equation}
\frac{\partial^2 \psi_k(r,z)}{\partial r^2}+\frac{1}{r}\frac
{\partial \psi_k(r,z)}{\partial r }- \frac{\psi_k(r,z )}{r^2} + 
\frac{1}{\lambda_k^2}\frac{\partial^2 \psi_k(r,z)}
{\partial z^2 }=0, \quad k=1,2,
\end{equation}
with pole $ (a,0) $ and are given by
\[ G_k(r,\lambda_k z; a,0)= \int\limits_0^\infty J_1(at)J_1(tr) e^{-t\lambda_k z}\,d\,t,  \quad  k=1,2. \]
Otherwise, they  can be rewritten as
\[ G_k(r,\lambda_k z; a,0)= \dfrac{1}{\pi}\int\limits_0^\pi \dfrac{\cos  \alpha }{(r^2+a^2+\lambda_k^2 z^2-2a r\cos \alpha )^{1/2}} \,d\,\alpha \]
The functions 
\begin{equation*}
F_k(r,z)= \int\limits_{0}^{\infty} e^{-t\lambda_k z}t J_1(at)J_0(tr)\,d\,t, k=1,2
\end{equation*}
by the help of the well-known relation 
\[ \frac{d}{d\zeta}(\zeta J_1(\zeta))= \zeta J_0(\zeta), \]
where $ \zeta=tr $, are expressed as
\begin{gather*}
 F_k(r,z)= \dfrac{1}{r}\int\limits_{0}^{\infty} e^{-t\lambda_k z}t J_1(at)J_0(tr)\,d\,t= \nonumber \\
\dfrac{1}{r}\int\limits_{0}^{\infty} e^{-t\lambda_k z}t J_1(at)\dfrac{d}{d(tr)}(tr)J_1(tr))\,d\,t= 
 \frac{1}{r}\frac{\partial}{\partial r}(r G_k(r,\lambda_k z)), k=1,2. 
\end{gather*} 
Displacements are convolutions
\begin{equation}
\begin{split}
 u_r(r,z)&= -\frac{k_2 \lambda_2}{\delta}(-\frac{\partial}{\lambda_1 \partial z})G_1(r,\lambda_1 z)\ast f_1 +  \frac{k_1 \lambda_1}{\delta}(-\frac{\partial}{\lambda_2 \partial z})G_2(r,\lambda_2 z)\ast f_1 \\
         & -\frac{1}{\delta}\frac{\partial}{\partial r}g_1(r,\lambda_1 z)\ast f_2 +  \frac{1}{\delta}\frac{\partial}{ \partial r}g_2(r,\lambda_2 z)\ast f_2,
\end{split}
\end{equation}
\begin{equation}
\begin{split}
u_z(r,z)&=-\frac{k_2 \lambda_2}{\delta}(-\frac{\partial}{\lambda_2 \partial z})g_2(r,\lambda_2 z)\ast f_2+
\frac{k_1 \lambda_1}{\delta}(-\frac{\partial}{\lambda_1 \partial z})g_1(r,\lambda_1 z)\ast f_2+\\
        &\frac{k_1k_2 \lambda_1 \lambda_2}{\delta}\{\frac{1}{r}\frac{\partial}{\partial r}(r G_2(r,\lambda_2 z))\ast f_1-\frac{1}{r}\frac{\partial}{\partial r}(r G_1(r,\lambda_1 z))\ast f_1\}.
\end{split}
\end{equation}
Here the $ \ast $-operation, called below a (Hankel) convolution, is defined as follows:
\[M(r,z)\ast f(r) = \int\limits_0^\infty M(r,z;a,0) f(a) a \,d\,a. \]

Now, the function
\[ -\dfrac{1}{\lambda_k} \dfrac{\partial}{\partial z_k} G_k(r,\lambda_k z;a,0)= \Lambda_{k1}(r, \lambda_k z;a,0) \]
is expressed as
\[ \Lambda_{k1}(r, \lambda_k z;a,0)= \dfrac{1}{\pi} \int\limits_0^\pi  \dfrac{\lambda_k z \cos \alpha }{(r^2 +a^2 +\lambda_k^2 z^2 -2 ra \cos \alpha)^{3/2}}\,d\,\alpha, \quad k=1,2.  \]
Then displacement $ u_r(r,z) $ is
\begin{gather}
u_r(r,z)= \Lambda_{11}(r,\lambda_1 z;a,0)*f_1(r) + \nonumber \\
\dfrac{k_1 \lambda_1}{\delta}\{ \Lambda_{21}(r,\lambda_2 z;a,0)*f_1(r)- \Lambda_{11}(r,\lambda_1 z;a,0)*f_1(r) \} + \nonumber \\
\frac{1}{\delta}\frac{\partial}{ \partial r}g_2(r,\lambda_2 z)\ast f_2(r)- \frac{1}{\delta}\frac{\partial}{\partial r}g_1(r,\lambda_1 z)\ast f_2(r).
\end{gather}
Here
\[ \dfrac{\partial}{ \partial r}g_k(r, \lambda_k z;a,0)=\dfrac{1}{\pi} \int\limits_0^\pi \dfrac{r-a\cos \alpha}{(r^2 +a^2 +\lambda_k^2 z^2 - 2ra \cos \alpha)^{3/2}} \,d\,\alpha, \quad k=1,2. \]
Integrals in (3.9) should be understood as singular principal value integrals. 
In a similar way,
\[ -\dfrac{\partial}{\lambda_k \partial z}g_k(r, \lambda_k z;a,0)=\dfrac{1}{\pi} \int\limits_0^\pi \dfrac{\lambda_k z}{(r^2 +a^2 +\lambda_k^2 z^2 - 2ra \cos \alpha)^{3/2}} \,d\,\alpha, \quad k=1,2. \]
and is the Poisson kernel for the equation (2.1). Displacement $ u_z(r,z) $ is
\begin{gather}
u_z(r,z)= -\dfrac{1}{\lambda_2}(\dfrac{\partial}{\partial z})g_2(r, \lambda_2 z)*f_2(r) + \nonumber \\
\dfrac{k_1 \lambda_2}{\delta}\{(-\dfrac{1}{\lambda_1}\dfrac{\partial}{\partial z}g_1(r, \lambda_1 z))*f_2(r)) -(-\dfrac{1}{\lambda_2}(\dfrac{\partial}{\partial z}g_2(r, \lambda_2 z; a,0))*f_2(r)) ) \} + \nonumber \\
\dfrac{k_1k_2 \lambda_1 \lambda_2}{\delta}\{ \dfrac{1}{r}\dfrac{\partial (r G_2(r, \lambda_2 z))}{\partial r}*f_1(r)- \dfrac{1}{r}\dfrac{\partial (r G_1(r, \lambda_1 z))}{\partial r}*f_1(r).  \}
\end{gather}
Convolutions with $ f_1(r) $ in (3.10) are singular integrals. It is obvious, nevertheless, that boundary conditions are satisfied, as differences of singular integrals vanish on the boundary. 

This is the formal solution of the problem (3.1)-(3.2). If boundary data are continuous, then displacements are continuous up to the boundary. To make its solution unique, one has to require boundedness of the strain energy; it means, that first order derivatives of displacements should be square integrable in the closed half-plane. Therefore, boundary data have to belong to
$ C^{0, \alpha}({\bf R}), \alpha > 1/2 $ on any compact subset of the real line, or have to satisfy the more stronger assumption: $ f_k(r) \in C^{1, \alpha}({\bf R}), k=1,2 $. The above told assumptions can be proved, but it is a theme for another paper. Here $ C^{m, \alpha}({\bf R}), m \geq 0 $ is the class of functions with $ m $ continuous derivatives with the $ m $-th order derivative satisfying the H\"{o}lder condition. Recall also, that the existence of a smooth solution of the Dirichlet problem for a single axially symmetric elliptic equation in a bounded region was proved in \cite{Par}. 

\section{The complex-valued presentation of the previous solution }
Introduce now the complex variable approach. Put
\[ \varphi_k(r,z) =\mathit{Re}\/ \Psi_k(r+i\lambda_k z), \qquad k=1,2, \]
then
\[ u_r(r,z)=\mathit{Re}\sum\limits_{j=1}^2 \Psi_j'(r+i\lambda_j z), \quad u_z(r,z) = \mathit{Re}\sum\limits_{j=1}^2i k_j \lambda_j \Psi_j'(r+i\lambda_k z). \]
Here $ ' $ means a derivative with an argument of a function.
Put also
\[ \Psi_1'(r+i\lambda_1 z)= m_1(r,\lambda_1 z)+i m_2(r,\lambda_1 z), \Psi_2'(r+i\lambda_2 z)= m_3(r,\lambda_2 z)+i m_4(r,\lambda_2 z). \]

The functions  $ m_k, k=1,2,3,4 $, are given by (for brevity, the dependence of fundamental solutions on pole $ (a,0) $ is suppressed)
\begin{equation*}
m_1(r,z)= -\dfrac{k_2\lambda_2}{\delta} (-\dfrac{1}{\lambda_1}\dfrac{\partial G_1(r,\lambda_1 z)}{\partial z})* f_1(r)- \dfrac{1}{\delta}\dfrac{\partial g_1(r,\lambda_1 z)}{\partial r}*f_2(r),
\end{equation*}
\begin{equation*}
m_3(r,z)= \dfrac{k_1\lambda_1}{\delta} (-\dfrac{1}{\lambda_2}\dfrac{\partial G_2(r,\lambda_2 z)}{\partial z})* f_1(r)+ \dfrac{1}{\delta}\dfrac{\partial g_2(r,\lambda_2 z)}{\partial r}*f_2(r),
\end{equation*}
\begin{equation*}
m_2(r,z)= -(-\dfrac{1}{\lambda_1}\dfrac{\partial g_1(r,\lambda_1 z)}{\partial z})*f_2(r)+ \dfrac{k_2 \lambda_2}{\delta}\dfrac{1}{r}\dfrac{\partial}{\partial r}(r G_1(r, \lambda_1 z))*f_1(r),
\end{equation*}

\begin{equation*}
m_4(r,z)= (-\dfrac{1}{\lambda_2}\dfrac{\partial g_2(r,\lambda_2 z)}{\partial z})*f_2(r)- \dfrac{k_1 \lambda_1}{\delta}\dfrac{1}{r}\dfrac{\partial}{\partial r}(r G_1(r, \lambda_2 z))*f_1(r).
\end{equation*}
As result, $ \tilde{P}_k(r+i\lambda_k z)= \Psi'_k(r+i \lambda_k z) $ are expressed as
\begin{gather*}
\tilde{P}_1(r+i\lambda_1 z)= -\dfrac{k_2\lambda_2}{\delta} (-\dfrac{1}{\lambda_1}\dfrac{\partial G_1(r,\lambda_1 z)}{\partial z})* f_1(r)- \dfrac{1}{\delta}\dfrac{\partial g_1(r,\lambda_1 z)}{\partial r}*f_2(r)+  \nonumber \\
i(\dfrac{k_1\lambda_1}{\delta} (-\dfrac{1}{\lambda_2}\dfrac{\partial G_2(r,\lambda_2 z)}{\partial z}* f_1(r)+ \dfrac{1}{\delta}\dfrac{\partial g_2(r,\lambda_2 z)}{\partial r}*f_2(r)),
\end{gather*}
\begin{gather*}
\tilde{P}_2(r+i\lambda_2 z)= \dfrac{k_1\lambda_1}{\delta} (-\dfrac{1}{\lambda_2}\dfrac{\partial G_2(r,\lambda_2 z)}{\partial z})* f_1(r)+ \dfrac{1}{\delta}\dfrac{\partial g_2(r,\lambda_2 z)}{\partial r}*f_2(r)+  \nonumber \\
i((-\dfrac{1}{\lambda_2}\dfrac{\partial g_2(r,\lambda_2 z)}{\partial z}*f_2(r)- \dfrac{k_1 \lambda_1}{\delta}\dfrac{1}{r}\dfrac{\partial}{\partial r}(r G_1(r, \lambda_2 z))*f_1(r)).
\end{gather*}
Put
\begin{gather}
N_k(r+ i\lambda_k z)= -\left( -\dfrac{1}{\lambda_k}\dfrac{\partial G_k(r, \lambda_k z)}{\partial z}\right)*\gamma_{k1} -\dfrac{\partial g_k(r,\lambda_k z)}{\partial r}*\gamma_{k2} +  \nonumber \\
i \left(-(-\dfrac{1}{\lambda_k}\dfrac{\partial g_k(r, \lambda_k z)}{\partial z})* \gamma_{k2} + \dfrac{1}{r}\dfrac{\partial (r G_k(r, \lambda_k z))}{\partial r}* \gamma_{k1}\right), \quad k=1,2,
\end{gather}
and solve the equations
\begin{equation*}
\mathit{Re}\sum_{j=1}^2 N_j(r+i \lambda_j z)\vert_{z=0}= f_1(r), \quad \mathit{Re}\sum_{j=1}^2i \lambda_j k_j N_j(r+i \lambda_j z)\vert_{z=0}= f_2(r),
\end{equation*}
with respect to $ \gamma_{ij}, i,j=1,2 $. Then
\begin{equation*}
\gamma_{11}= -\dfrac{k_2 \lambda_2}{\delta}f_1(r), \gamma_{21}= \dfrac{k_1 \lambda_1}{\delta}f_1(r), \gamma_{12}= -\dfrac{1}{\delta}f_2(r), \gamma_{22}= \dfrac{1}{\delta}f_2(r),
\end{equation*}
i.e., we derive the previous formulas for $ \tilde{P}_k(r+ i\lambda_k z), k=1,2. $ The functions
$ N_k(r+i \lambda_k z), k=1,2 $ are the  required axially symmetric Cauchy type integrals for an upper half-plane. Now, it is necessary to rewrite these formulas for an arbitrary bounded axially symmetric region. This is done in the next section.

\section{Axially symmetric potentials for a bounded region}
\subsection{The Cauchy type integral for a single axially symmetric equation}

Let  $ \Phi(r,\phi,z) $ be a harmonic function. Written in the cylindrical coordinates, it is a solution of the equation
\[ \dfrac{\partial^2 \Phi}{\partial r^2}+ \dfrac{1}{r}\dfrac{\partial \Phi}{\partial r}+ \dfrac{1}{r^2}\dfrac{\partial^2 \Phi}{\partial \phi^2}+ \dfrac{\partial^2 \Phi}{\partial z^2}=0. \]
Expand  $ \Phi(r,\phi,z) $ in a Fourier series with respect to $ \phi $, that is, write it as a series
\[ \Phi(r,\phi,z)= \Phi_0(r,z) + \sum\limits_{n=1}^\infty \{\Phi_n^c(r,z)\cos n \phi + \Phi_n^s \sin n \phi \}.  \]
Then functions $ \Phi_n^c, \Phi_n^s $ are solutions of the equation
\begin{equation}
\dfrac{\partial^2 \Phi_n}{\partial r^2}+ \dfrac{1}{r}\dfrac{\partial \Phi_n}{\partial r} - \dfrac{n^2}{r^2} \Phi_n + \dfrac{\partial^2 \Phi_n}{\partial z^2}=0, n=0,1, \ldots.
\end{equation}
In particular, 
\begin{equation}
\dfrac{\partial^2 \Phi_0}{\partial r^2}+ \dfrac{1}{r}\dfrac{\partial \Phi_0}{\partial r} + \dfrac{\partial^2 \Phi_0}{\partial z^2}=0, 
\end{equation}
that is, the zero order term  $ \Phi_0(r,z) $  is an axially symmetric harmonic function. Now, if  $ R $ is the distance between points with Cartesian coordinates $ (x_1,x_2,x_3) $ and $ (x_{10}, x_{20},x_{30}) $, where 
\[ x_1=r \cos \theta, \quad x_2= r\sin \theta, \quad x_3=z, \]
\[ x_{10}=r_1\cos\theta_0, \quad x_{20}=r_1 \sin\theta_0, \quad x_{30} =z_1, \]
then 
\[ R = \sqrt{r^2+r_1^2 +(z-z_1)^2 -2r r_1 \cos(\theta - \theta_0)}. \]
Put $ \phi= \theta- \theta_0 $, and 
expand  $ R^{-1} $ (the fundamental solution of the Laplace's equation) into the Fourier series with respect to the angular coordinate $ \phi $. Then  
\[ \dfrac{1}{R}= \dfrac{1}{\sqrt{2 r r_1}} \dfrac{1}{\sqrt{\cosh \psi -\cos  \phi}}= \sum\limits_{n=0}^\infty w_n(r,z)\cos n \phi.\]
Here
\[ 2r r_1 \cosh \psi = r^2 +r_1^2 +(z -z_1)^2. \]

The functions $ w_n(r,z;r',z'), n \geq 0 $ are the Fourier coefficients in the previous expansion, and can be written by  evenness of the integrand with respect to $  \phi $ as
\[ w_0 (r,z;r',z')= \dfrac{1}{ \pi}\int\limits_{0}^{\pi} \dfrac{\,d \,\phi}{\sqrt{r^2 + r_1^2 -2 r r_1 \cos \phi +(z-z_1)^2}}, \]
\[ w_n (r,z;r',z')=\dfrac{1}{ \pi}\int\limits_{0}^\pi \dfrac{  \cos n\phi \,d\,\varphi}{\sqrt{r^2 + r_1^2 -2 r r_1 \cos \phi +(z-z_1)^2}}, \quad n \geq 1. \]
Therefore, $ w_n(r,z; r',z'), n \geq 0 $ is the axially symmetric fundamental solution of the equation (5.1) with pole at $ (r',z') $, corresponding to the $n$-th harmonic of this expansion. Now, D. I. Sherman
represented a harmonic function $ u(x_1,x_2, x_3) $ as a simple layer potential 
\[ u(x_1,x_2,x_3)= \dfrac{1}{2 \pi}\int\limits_{\partial Q} \nu(M_1) \dfrac{\,d\,\sigma_1}{\varrho_{10}},  \]
where $ \varrho_{10} $ is the distance between points $ M_0(r, \phi, z) $ and $ M_1(r_1, \phi_1, z_1) $, and $ r_1=r(s), \phi_1, z_1=z(s) $ are variables of integration. If density $ \nu(M_1) $ is independent of the angular coordinate $ \phi $, then
\begin{gather*}
\dfrac{1}{2 \pi} \int\limits_{\partial Q} \nu(M_1) \dfrac{\,d\,\sigma_1}{\varrho_{10}}= \dfrac{1}{2 \pi}\int\limits_{L} \nu(r(s),z(s)) \\
\int\limits_0^{2\pi}\dfrac{\,d\,\phi_1}{\sqrt{(z_1(s)-z)^2 +r_1(s)^2 +r^2 -2r r_1(s) \cos (\phi_1  -\phi)}} r(s) \,d\,s.
\end{gather*}
Here $ s $ is the arc length parameter of the meridional section of a solid, $ \,d\,s $ is the arc length element, $ L $ is the meridional section of a boundary, given by equations $ r=r(s), z=z(s) $. The internal integral is independent of $\phi $, and is equal to \newline
$ w_0(r,z; r(s),z(s)) $. Therefore, we arrive to the simple layer axially symmetric potential. In a similar way can be written the double layer potential. In short, it can be written as
\begin{equation}
u(x_1,x_2,x_3)= \dfrac{1}{2 \pi}\int\limits_{\partial Q} \nu(M_1) \dfrac{\,d\,\sigma_1}{\varrho_{10}}=\dfrac{1}{\pi}\int\limits_{L_1} \nu(s) w_0(r-r(s), z-z(s))\,d S,
\end{equation}
where $ \,dS= r(s) \,d\,s $ is the superficial axially symmetric area element. In a similar way can be written a double layer axially symmetric potential. Indeed, let $ {\bf n}=(-z'(s), r'(s)) $ be an internal normal vector. Then the normal derivative of a function $ v(r,z) $ is
\[ \dfrac{\partial v(r,z)}{\partial n}= \dfrac{\partial v}{\partial r} n_1 + \dfrac{\partial v}{\partial z}n_2. \]
Put $ v =w_0(r-r(s), z-z(s)) $. Then a double layer potential  $ U_0(r,z) $ is written as
\[ U_0(r,z) = \int\limits_{\partial Q}\mu(s) \dfrac{\partial w_0(r-r(s), z-z(s))}{\partial n} \,d\,S, \]
where $ \mu(s) $ is some unknown density. Indeed, it is proved in \cite{Raja}[Theorem 2], that for an axially symmetric region with the Lyapunov boundary when a point $ (r,z) $ tends to a point of its surface $( r(s_0),z(s_0)) $ from inside a region, the double layer potential has the value
\begin{equation}
U_0^+(r(s_0),z(s_0))= \mu(s_0) + U_0(r(s_0),z(s_0))
\end{equation}
in the limit, if its density $ \mu(s) $ is a continuous function.

In a similar way can be computed the tangential derivative of a simple layer potential. Recall, that the tangential derivative of a function $ v(r,z) $ is defined as
\[ \dfrac{\partial v(r,z)}{\partial \tau}= \dfrac{\partial v}{\partial r} \tau_1 + \dfrac{\partial v}{\partial z}\tau_2, \]
where $ \tau_1=r'(s), \tau_2= z'(s) $. Then the tangential derivative of the fundamental solution with pole at the point  $ (r(s),z(s)) $ is given by
\[ U_{00} (r,z) = \int\limits_{\partial Q}\mu(s) \dfrac{\partial w_0(r-r(s), z-z(s))}{\partial \tau} \,d\,S. \] 
In distinction to a double layer potential, a tangential derivative of a simple layer is a continuous function when crossing a boundary.
Now, form the sum
\begin{gather}
V(r,z)= U_0(r,z) +i U_{00}(r,z)= \nonumber \\
\int\limits_{\partial Q}\mu(s) \left( \dfrac{\partial w_0(r,r(s), z-z(s))}{\partial n} + i \dfrac{\partial w_0(r,r(s), z-z(s))}{\partial s}\right) \,d\,S.
\end{gather}
Then 
\begin{equation}
V(r(s_0),z(s_0))^+ =  \mu(s_0) + V(r(s_0),z(s_0)),
\end{equation}
when a point $ (r,z) $ tends to a point of its surface $( r(s_0),z(s_0)) $ from inside the region. Here the superscript $ + $ denotes the limit value of a function at a boundary. But it is not the whole truth, as we have to deal with an axially symmetric Cauchy-Riemann system of functions 
\begin{equation}
\frac{1}{r}\frac{\partial q_1}{\partial z}-\frac{\partial q_2}{\partial r}=0,
\end{equation}
\begin{equation}
 \frac{1}{r}\frac{\partial q_1}{\partial r}+\frac{\partial q_2}{\partial z}=0.
\end{equation}
Here $ q_1,q_2 $ are solutions of two different elliptic equations:
\begin{equation}
 \dfrac{1}{\partial r}\left(\dfrac{1}{r}\dfrac{\partial q_1}{\partial r} \right) + \dfrac{1}{r}\dfrac{\partial^2 q_1}{\partial z^2}=0,
\end{equation}
\begin{equation}
\dfrac{\partial}{\partial r}\left( r\dfrac{\partial q_2}{\partial r}\right)+ r \dfrac{\partial^2 q_2}{\partial z^2}=0.
\end{equation}
The systems of equations, more general, than the system (5.7), (5.9) was studied by A. Weinstein in \cite{Wein} and his numerous papers under the name \textquotedblleft Generalized axially symmetric potential theory\textquotedblright. 
The equation (5.10) was considered above. Consider  now the equation (5.9). A function $ q_1(r,z) $ can be written as a product $ q_1=rq_3 $, where $ q_3 $ is a solution of the equation
\begin{equation}
\dfrac{\partial^2 q_3(r,z)}{\partial z^2}+ \dfrac{\partial^2 q_3(r,z)}{\partial r^2}+ \dfrac{1}{r}\dfrac{\partial q_3(r,z )}{\partial r}- \dfrac{q_3(r,z)}{r^2}=0.
\end{equation}
The function 
\[ w_1 (r,z;r(s),z(s))=\dfrac{1}{ \pi}\int\limits_{0}^\pi \dfrac{  \cos \phi \,d\,\varphi}{\sqrt{r^2 + r_1^2 -2 r r_1 \cos \phi +(z-z_1)^2}}  \]
is the axially symmetric fundamental solution with pole at the point $ (r_1,z_1) $, where $r_1= r(s),z_1= z(s) $, of this equation. Compute the functions
\begin{equation}
\dfrac{1}{r}\dfrac{\partial}{\partial z}(r w_1(r,z;r_1(s),z(s)), \quad \dfrac{1}{r}\dfrac{\partial}{\partial r}(r w_1(r,z;r(s),z(s)) 
\end{equation}
and construct the corresponding normal and tangential derivatives, that is,
\begin{gather}
 \dfrac{1}{r}\dfrac{\partial}{\partial n}(r w_1(r,z;r(s),z(s))= \nonumber \\
\dfrac{1}{r}\dfrac{\partial}{\partial r}(r w_1(r,z;r(s),z(s))n_1 +
\dfrac{1}{r}\dfrac{\partial}{\partial z}(r w_1(r,z;r(s),z(s))n_2, 
\end{gather}
\begin{gather}
\dfrac{1}{r}\dfrac{\partial}{\partial \tau}(r w_1(r,z;r(s),z(s))= \nonumber \\
\dfrac{1}{r}\dfrac{\partial}{\partial r}(r w_1(r,z;r(s),z(s))\tau_1 +
\dfrac{1}{r}\dfrac{\partial}{\partial z}(r w_1(r,z;r(s),z(s))\tau_2.
\end{gather}
and form the sum
\begin{gather}
V_{00}(r,z)= \nonumber \\
\int\limits_{\partial Q}\mu_1(s)\left( \dfrac{1}{r}\dfrac{\partial}{\partial n}(r w_1(r,z;r(s),z(s))+  
i \dfrac{1}{r}\dfrac{\partial}{\partial \tau}(r w_1(r,z;r_1(s),z_1(s))) \right) \,d\,S.
\end{gather}
Introduce for brevity shorter notation, namely, put
\begin{equation}
G(r,z,s)= w_1(r,z;r(s),z(s))) , g(r,z,s)= w_0(r,z;r(s),z(s)).
\end{equation}
Now, form the sum
\begin{gather}
\Pi(r,z)=V_{00}(r,z)+ i U_{00}(r,z)= \nonumber \\
\int\limits_{\partial Q}\mu_1(s)\left( \dfrac{1}{r}\dfrac{\partial}{\partial n}(r w_1(r,z;r(s),z(s))+  
i \dfrac{1}{r}\dfrac{\partial}{\partial \tau}(r w_1(r,z;r(s),z(s))) \right) \,d\,S + \nonumber \\
i \int\limits_{\partial Q}\mu(s) \left( \dfrac{\partial w_0(r,r(s), z-z(s))}{\partial n} + i \dfrac{\partial w_0(r,r(s), z-z(s))}{\partial \tau}\right) \,d\,S.
\end{gather} 
It is the required Cauchy type integral. It is obvious, that the limit of the function $ \Pi(r,z) $, when a point $ (r,z) $ tends to a point $ (r(s_0),z(s_o)) $ inside a region, is equal to
\begin{equation}
\Pi^+ (r(s_0), z(s_0)) =\mu_1(s_0)+i \mu(s_0) + \Pi(r(s_0), z(s_0)), 
\end{equation}
if functions $ \mu(s) , \mu(s) $ are H\"{o}lder continuous.
This integral is similar to the standard Cauchy type integral in the usual complex variable and has, in general, the same jump properties when crossing the boundary. The formula (5.18) can be written in the other way. We have
\begin{gather}
\int\limits_{\partial Q}\mu_1(s)\left( \dfrac{1}{r}\dfrac{\partial}{\partial n}(r w_1(r,z;r(s),z(s))+  
i \dfrac{1}{r}\dfrac{\partial}{\partial \tau}(r w_1(r,z;r(s),z(s))) \right) \,d\,S = \nonumber \\
i\int\limits_{\partial Q}\mu_1(s)\left( \dfrac{1}{r}\dfrac{\partial}{\partial r}(r w_1(r,z;r(s),z(s))-  
i \dfrac{1}{r}\dfrac{\partial}{\partial z}(r w_1(r,z;r(s),z(s))) \right) r(s) \,d\,t, \nonumber \\
t(s)=r(s)+iz(s), \quad    t'(s)=(r'(s)+i z'(s))\,d\,s.
\end{gather}
In a similar way,
\begin{gather}
i \int\limits_{\partial Q}\mu(s) \left( \dfrac{\partial w_0(r,r(s), z-z(s))}{\partial n} + i \dfrac{\partial w_0(r,r(s), z-z(s))}{\partial \tau}\right) \,d\,S = \nonumber \\
-\int\limits_{\partial Q}\mu(s) \left( \dfrac{\partial w_0(r,r(s), z-z(s))}{\partial r} - i \dfrac{\partial w_0(r,r(s), z-z(s))}{\partial z}\right)r(s) \,d\,t.
\end{gather}
Introduce now the functions
\begin{gather*}
v_1(r,z;r(s),z(s))=\int\limits_{\partial Q}( \gamma_1 \{\dfrac{1}{r}\dfrac{\partial (r G(r,z;r(s),z(s)))}{\partial z}n_2 +\dfrac{1}{r}\dfrac{\partial rG(r,z;r(s), z(s))}{\partial r}n_1 \}+    \\ 
\gamma_2\{ \dfrac{1}{r}\dfrac{\partial (rG(r,z;r(s),z(s))}{\partial z}n_1 -\dfrac{1}{r}\dfrac{\partial (rG(r,z;r(s),z(s))}{\partial r}n_2 \}r(s) \,d\,s,
\end{gather*}
\begin{gather*}
 v_2(r,z)= \int\limits_{\partial Q} \Big(\gamma_1\lbrace \dfrac{\partial g(r,z;r(s),z(s))}{\partial r}n_2 -\dfrac{\partial g(r,z;r(s),z(s))}{\partial z}n_1\rbrace+ \\
 \gamma_2 \lbrace \dfrac{\partial G(r,z;r(s),z(s)}{\partial r}n_1 + \dfrac{\partial g(r,z;r(s),z(s) )}{\partial z}n_2\rbrace \Big)r(s)\,d\,s,
\end{gather*}
and take the  sum $ M(r,z)=v_1(r,z)+iv_2(r,z), i=\sqrt{-1} $. The result of summation is written as
\begin{gather*}
M(r,z; r(s), z(s))= \\
\int\limits_{\partial Q}\left( \gamma_1 \dfrac{1}{r}\dfrac{\partial (r G(r,z;r(s),z(s))}{\partial n} +   
\gamma_2 \dfrac{1}{r}\dfrac{\partial (rG(r,z;r(s),z(s)))}{\partial \tau}\right) r(s) \,d\,s + \nonumber \\
i\int\limits_{\partial Q} \left(\gamma_1 \dfrac{\partial g(r,z;r(s),z(s)}{\partial \tau} +\gamma_2 \dfrac{\partial g(r,z;r(s, z(s)))}{\partial n}\right) r(s) \,d\,s.
\end{gather*}
Denote the limit of the function $ M(r,;r(s),z(s)) $ when it tends to a boundary point from inside a region, by $ M^+(r(s_0),z(s_0)) $. Then
\[  M^+(r(s_0),z(s_0))= \gamma_1(s_0)+ i \gamma_2(s_0) + M(r(s_0), z_0) \]
on $ \partial Q $. Nevertheless, it is not the required Cauchy type integral. Recall, that for the equation 
\begin{equation*}
\dfrac{\partial^2 U_\nu}{\partial r^2}+ \dfrac{ \nu}{r}\dfrac{\partial U_\nu}{\partial r}+\dfrac{\partial^2 U_\nu}{\partial z^2}=0,  \quad \nu \geq 0 
\end{equation*}
the complete potential theory in H\"{o}lder classes of functions was worked out in \cite{Raja}. Some results from that  paper were used above. In particular, the Dirichlet and Neumann problems was solved in bounded axially symmetric regions for the previous equation.

\subsection{Axially Symmetric Potentials in transformed regions  }
Let $ Q $ be a bounded open symmetric region in the $ (r,z) $ plane which is symmetric with respect to the $ z $ axis and is such that its intersection with the $ z $ axis is an open interval. Let $ Q_{+} = Q \cap \{(r,z), r > 0 \} $. Let $ \Gamma_{+} $ denote the boundary of $ Q_+ $, and let $ \Gamma_{+} =\Gamma_1 \cup \Gamma_2 $, where $ \Gamma_1 $ is an interval on the axis and $ \Gamma_2 $ is a smooth Lyapounov curve lying entirely in the half-plane $ r > 0 $. Introduce the coordinates $ (r, \lambda_j z),j=1,2 $. It is obvious, that to any point $ (r,z) \in  Q $ corresponds a unique point $ (r, \lambda_j z), j=1,2 $, belonging to the region $ Q_j = \{(r,z_j), z_j=\lambda_j z, j=1,2 \} $ if $ (r,z) \in Q $.
Regions $ Q_j, j=1,2 $, following the analogy with the two-dimensional anisotropic theory of elasticity, can be called the affinely-transformed ones. Obviously, their topological and differential properties are identical.
Let $ g(r,z) \in C( \bar Q) $. Then there exists a unique axially symmetric potential $ \varphi(r,z)\in C^2(\bar Q)\cap C( \bar Q) $ which equals to $ g(r,z) \in C( \bar Q) $ on $ \Gamma_2 $. Moreover, it is a real analytic function for all interior points of $ Q $.
We write now the required potentials  in the regions $ Q_j, j=1,2 $.

Therefore, $ g_k(r, \lambda_k z;a,0)  = w_o(r,\lambda_k z;a,0)$, $ G_k(r, \lambda_k z;a,0)  = w_1(r,\lambda_k z;a,0)$. Now, if a meridional section of a region $ Q $  is described by the equation
$ ds =\sqrt{dr^2+dz^2} $, then the meridional section of the region $ Q_k $ is described by the equation $ ds_k= \sqrt{dr^2 + \lambda_k^2 dz^2} $, or, in the complex form,
\[ dt=(r'(s)+iz'(s))\,d\,s, \quad dt_k=(r'(s)+i\lambda_k z'(s))\,d\,s , \quad k=1,2. \]
Put, for brevity, 
\[ g_k(r, \lambda_k z)=w_o(r,\lambda_k z;r(s),\lambda_k z(s)), \quad G_k(r, \lambda_k z)=w_1(r,\lambda_k z;r(s),\lambda_k z(s)), \]
and compute the corresponding normal and tangential derivatives in the regions $ Q_k, k=1,2 $.

Here
\begin{gather}
 \dfrac{1}{r}\dfrac{\partial(r G_k(r, \lambda_k z))}{\partial n}=  \dfrac{1}{r}\dfrac{\partial(r G_k(r,\lambda_k z))}{\lambda_k \partial z}n_{2k} +\dfrac{1}{r}\dfrac{\partial(r G_k(r,\lambda_k z))}{\partial r}n_{1k}, \\
 \dfrac{1}{r}\dfrac{\partial(r G_k(r, \lambda_k z))}{\partial \tau}= \dfrac{1}{r}\dfrac{\partial(r G_k(r,\lambda_k z))}{\lambda_k \partial z}n_{1k} -\dfrac{1}{r}\dfrac{\partial(r G_k(r,\lambda_k z))}{\partial r}n_{2k}, \\
 \dfrac{\partial g_k(r,\lambda_k z)}{\partial \tau}= - \dfrac{\partial g_k(r,\lambda_k z) }{\partial r}n_{2k} +\dfrac{\partial g_k(r,\lambda_k z)}{\lambda_k \partial z}n_{1k} , \\
 \dfrac{\partial g_k(r,\lambda_k z)}{\partial n} =  \dfrac{\partial g_k(r, \lambda_k z)}{\partial r}n_{1k} + \dfrac{\partial g_k(r, \lambda_k z)}{\lambda_k \partial z}n_{2k}, \quad  k=1,2,  
\end{gather}
where 
\[ n_{1k}=\dfrac{-\lambda_k z'(s)}{\sqrt{r'(s)^2 + \lambda_k^2 z'(s)^2}}, \quad n_{1k}=\dfrac{r'(s)}{\sqrt{r'(s)^2 + \lambda_k^2 z'(s)^2}}, \quad k=1,2, \]
and are normal direction cosines of a boundary in the coordinate system $ (r,\lambda_k z) $. Therefore, the required Cauchy type integrals in the transformed regions are given by
\begin{gather}
\Pi_k(r+ i \lambda_k z)= \nonumber \\
\int \limits_{\partial Q} ( \dfrac{1}{r}\dfrac{\partial(r G_k)(r,\lambda_k z;r_1(s)),z_1(s)}{\partial n}
\gamma_{k1}(s)- \dfrac{\partial g_j(r,\lambda_j z; r(s), z(s))}{\partial \tau}\gamma_{k2}(s))r(s)\,d\,S_k +  \nonumber \\
 i\int\limits_{\partial Q} (\dfrac{\partial g_k(r, \lambda_k z; r(s), z(s))}{\partial n}\gamma_{k2}(s)+ \dfrac{1}{r}\dfrac{\partial (r G_k)(r,\lambda_k z; r(s), z(s))}{\partial \tau}\gamma_{j1}(s))r(s) \,d\,S_k, k=1,2,
\end{gather}
where $ r(s) \,d\,S_k= r(s) \sqrt{ r'^2(s) + \lambda_k^2 z'(s)}ds, k=1,2 $, and where the normal and tangential derivatives are computed according to the formulas (5.23)--(5.26). There is one more equivalent form of the Cauchy type integral:
\begin{gather}
\Pi_k(r,\lambda_k z)=\int \limits_{\partial Q} ( \dfrac{1}{r}\dfrac{ \partial (r G(r,z))}{\partial z}
 + i \dfrac{1}{r}\dfrac{ (\partial r G(r,z))}{\partial r})\gamma_{1k} r(s)\,d\,t_k + \nonumber \\ 
 i \int\limits_{\partial Q} (\dfrac{\partial g(r,\lambda_k z)}{\partial z}+ i \dfrac{\partial g(r,z)}{\partial r})\gamma_{2k}r(s)\,d\,t_k, \quad k=1,2.
\end{gather}
Here $ t_k(s)= r(s)+i\lambda_k z(s)$, $ d t_k(s)= (r'(s)+i\lambda_k z'(s)) ds, k=1,2 $ and $\gamma_{jk}(r), j,k=1,2$ are unknown real H\"{o}lder continuous densities.

\section{The regular system of integral equations}
Applying the approach, developed by us in Section 4, we seek the functions $ \varphi_k(r,z), k=1,2 $ and displacements $ u_r(r,z), u_z(r,z) $ as
\[ \varphi_k(r,z) =\mathit{Re}\/ \Psi_k(r+i\lambda_k z), \qquad k=1,2, \]
then
\[ u_r(r,z)=\mathit{Re}\sum\limits_{j=1}^2 \Psi_j'(r+i\lambda_j z), \quad u_z(r,z) = \mathit{Re}\sum\limits_{j=1}^2i k_j \lambda_j \Psi_j'(r+i\lambda_k z). \]
Put
\[ \Psi_k'(r+i\lambda_z )= \Pi_k(r, \lambda_k z;r(s), \lambda_k z(s)), \quad k=1,2,  \]
where $ \Pi_k(r, \lambda_k z;r(s), \lambda_k z(s)), k=1,2 $ are determined by formulas (5.27).

Now we can derive the required system of integral equations in a bounded axially symmetric region with a smooth boundary.
Determine real functions $ \gamma_{ij}, i,j=1,2 $, entering into (5.27), from  the equations
\begin{gather}
\gamma_{11}+\gamma_{21}= h_1(s), \quad k_1\lambda_1 \gamma_{11}+ k_2\lambda_2 \gamma_{21}=0, \nonumber \\
\gamma_{12}+\gamma_{22}=0, \quad -k_1 \lambda_1 \gamma_{12} -k_2 \lambda_2 \gamma_{22}= h_2(s).
\end{gather}
As result, the displacements are
\begin{gather*}
u_r(r,z)= -\dfrac{k_2 \lambda_2}{\delta} \int\limits_{\partial Q} h_1(s) \dfrac{1}{r}\dfrac{\partial(r G_1(r, \lambda_1 z))}{\partial n}r(s)\,d\,s + \nonumber \\
\dfrac{k_1 \lambda_1}{\delta} \int\limits_{\partial Q} h_1(s)\dfrac{1}{r}\dfrac{\partial(r G_2(r, \lambda_2 z))}{\partial n}r(s)\,d\,s +  \nonumber \\
\dfrac{1}{\delta}\int\limits_{\partial Q}h_2(s) (\dfrac{\partial g_1(r, \lambda_1 z)}{\partial \tau}-\dfrac{\partial g_2(r, \lambda_2 z)}{\partial \tau })r(s)\,d\,s,
\end{gather*}
\begin{gather*}
u_z(r,z)= \dfrac{k_1 \lambda_1}{\delta}\int\limits_{\partial Q} h_2(s) \dfrac{\partial(r g_1(r,\lambda_1 z))}{\partial n}r(s)\,d\,s- \dfrac{k_2 \lambda_2}{\delta}\int\limits_{\partial Q} h_2(s) \dfrac{\partial g_2(r,\lambda_2 z)}{\partial n}r(s)\,d\,s+  \nonumber \\
\dfrac{k_1 k_2 \lambda_1 \lambda_2}{\delta}\int\limits_{\partial Q} h_1(s)(\dfrac{1}{r}\dfrac{\partial(r G_1)(r, \lambda_1 z))}{\partial \tau}- \dfrac{\partial (r G_2(r, \lambda_2 z))}{\partial \tau})r(s)\,d\,s.
\end{gather*}
Recall, that here the functions $ g_k(r,\lambda_k z; r(s), \lambda_k z(s)), G_k(r,\lambda_k z; r(s), \lambda_k z(s)) $ depend on $ r(s), z(s) $, but for brevity this dependence is not indicated explicitly. The previous formulas can be rewritten as
\begin{gather}
u_r(r,z)= \int\limits_{\partial Q} h_1(s)\dfrac{\partial( r G_1(r, \lambda_1 z))}{\partial n}r(s)\,d\,s+ \nonumber \\
\dfrac{k_1 \lambda_1}{\delta} \int\limits_{\partial Q} h_1(s)\big( \dfrac{\partial(r G_2(r, \lambda_2 z))}{\partial n} -\dfrac{\partial(r G_1(r, \lambda_1 z))}{\partial n} \big)r(s)\,d\,s+  \nonumber \\
\dfrac{1}{\delta}\int\limits_{\partial Q}h_2(s) (\dfrac{\partial g_1(r, \lambda_1 z)}{\partial \tau}-\dfrac{\partial g_2(r, \lambda_2 z}{\partial \tau})r(s)\,d\,s,
\end{gather}
\begin{gather}
u_z(r,z)= \int\limits_{\partial Q} h_2(s) \dfrac{\partial g_2(r,\lambda_2 z)}{\partial n}r(s)\,d\,s+ \nonumber \\
 \dfrac{k_1 \lambda_1}{\delta}\int\limits_{\partial Q} h_2(s) \big(\dfrac{\partial G_1(r,\lambda_1 z)}{\partial n} -  \dfrac{\partial g_2(r,\lambda_2 z)}{\partial n}\big)r(s) \,d\,s+  \nonumber \\
\dfrac{k_1 k_2 \lambda_1 \lambda_2}{\delta}\int\limits_{\partial Q} h_1(s)(\dfrac{1}{r}\dfrac{\partial(r G_1(r, \lambda_1 z))}{\partial \tau}- \dfrac{1}{r}\dfrac{\partial(r G_2(r, \lambda_2 z))}{\partial \tau})r(s)\,d\,s.
\end{gather}
Thus, densities $ h_k(s), k=1,2 $, are solutions of two integral equations
\begin{gather}
h_1(s_0)+ \int\limits_{\partial Q} h_1(s)\dfrac{1}{r}\dfrac{\partial(r G_1(r, \lambda_1 z))}{\partial n^0}r(s)\,d\,s + \nonumber \\
\dfrac{k_1 \lambda_1}{\delta} \int\limits_{\partial Q} h_1(s)\big( \dfrac{1}{r}\dfrac{\partial(r G_2(r, \lambda_2 z))}{\partial n^0} -\dfrac{1}{r}\dfrac{\partial(r G_1(r, \lambda_1 z))}{\partial n^0} \big)r(s)\,d\,s+  \nonumber \\
\dfrac{1}{\delta}\int\limits_{\partial Q}h_2(s) (\dfrac{\partial g_1(r, \lambda_1 z)}{\partial \tau^0}-\dfrac{\partial g_2(r, \lambda_2 z}{\partial \tau^0})r(s)\,d\,s = g_1(s_0),
\end{gather}
\begin{gather}
h_2(s_0)+ \int\limits_{\partial Q} h_2(s) \dfrac{\partial g_2(r,\lambda_2 z)}{\partial n^0}r(s)\,d\,s+ \nonumber \\
 \dfrac{k_1 \lambda_1}{\delta}\int\limits_{\partial Q} h_2(s) \big(\dfrac{g_1(r,\lambda_1 z)}{\partial n^0} -  \dfrac{\partial g_2(r,\lambda_2 z)}{\partial n^0}\big)r(s) \,d\,s+  \nonumber \\
\dfrac{k_1 k_2 \lambda_1 \lambda_2}{\delta}\int\limits_{\partial Q} h_1(s)(\dfrac{1}{r}\dfrac{\partial(r G_1(r, \lambda_1 z))}{\partial \tau^0}- \dfrac{1}{r}\dfrac{\partial(r G_2(r, \lambda_2 z)}{\partial \tau^0})r(s)\,d\,s= g_2(s_0).
\end{gather}
The superscript $ 0 $ above means the direct value of a potential at a boundary $ \partial Q $. It is the correct set of integral equations for the displacement problem, in distinction to the set, given in [ \cite{Alex}, p. 339, the formula (40.27)]. Obviously, formulas (6.2), (6.3) give the correct solution of the considered problem in the half-plane $ R_+^2=\{(r,z), z>0 \} $ (compare the result with the formula (4.1) for $ N_k(r+i\lambda_k z), k=1,2 $)
and give the system of regular integral equations for axially symmetric regions with a smooth boundary. The regularity of the system (6.4), (6.5)) can be proved similar to the approach, used  earlier in \cite{Bog2}, \cite{Bog3}.

\end{document}